\RequirePackage{fix-cm}
\documentclass[smallextended]{svjour3}       
\smartqed  
\usepackage{graphicx}
\usepackage[all, knot]{xy}
\usepackage{amsfonts,amssymb}
\usepackage{latexsym}
\usepackage{amsmath}
\usepackage[all]{xy}
\usepackage{stmaryrd}
\newtheorem{axiom}{Axiom}
\usepackage{url}
\begin{document}

\title{How Logic Interacts with Geometry: Infinitesimal Curvature of 
Categorical Spaces} 


\author{Michael Heller         \and
        Jerzy Kr\'ol 
}


\institute{M. Heller \at
              Copernicus Center for Interdisciplinary Studies, ul. S{\l}awkowska 17, 31-016 Krak\'ow, Poland \\
              \email{mheller@wsd.tarnow.pl}           
           \and
           J. Kr{\'o}l \at
              Institute of Physics, University of Silesia, ul. Uniwersytecka 
4, 40-007 Katowice, Poland\\
              \email{jerzy.krol@us.edu.pl}
}

\date{Received: date / Accepted: date}

\maketitle

\begin{abstract}
 In category theory, logic and geometry cooperate with each other producing 
what is known under the name Synthetic Differential Geometry (SDG). The main 
difference between SDG and standard differential geometry is that the  
intuitionistic logic of SDG enforces the existence of infinitesimal objects 
which essentially modify the local structure of spaces considered in SDG. We 
focus on an ``infinitesimal version'' of SDG, an infinitesimal 
$n$-dimensional formal manifold, and develop differential geometry on it. In 
particular, we show that the Riemann curvature tensor on infinitesimal level is itself infinitesimal. We construct a heuristic model $S^3 \times 
\mathbb{R} \subset \mathbb{R}^4$ and study it from two perspectives: the 
perspective of the category SET and that of the so-called topos $\mathcal{G}$ 
of germ-determined ideals. We show that the fact that in this model the curvature tensor 
is infinitesimal (in $\mathcal{G}$-perspective) eliminates the existing 
singularity. A surprising effect is that the hybrid geometry based on the existence of the infinitesimal 
and the SET levels generates an exotic smooth structure on $\mathbb{R}^4$. We briefly 
discuss the obtained results and indicate their possible applications.
\end{abstract}

\section{Introduction}
\label{intro}
Classical logic is a formalization of our everyday patterns of reasoning. In 
doing science and trying to understand the universe, we extend this way of 
reasoning to the entire realm of reality. The first serious warning that this 
huge extrapolation can be misleading came from quantum mechanics in which some 
time honoured principles of classical logic turned out to be invalid. An 
attempt to cope with this ``deviation'' led to what one calls quantum logic 
(see, for instance, \cite{Cohen,vNeumann,Omnes}). The next step was to change 
the ``deviation'' into a rule, and to place quantum mechanics in a topos 
environment where the modification of logic is an element of the game. The 
original Isham and Butterfield's idea \cite{IshamButter98,IshamButter99} 
has developed into a rich program of reconceptualization of quantum mechanics in 
terms of topos theory \cite{DoringIsham,Flori,Isham}. In this conceptual 
setting, intuitionistic logic provides a natural way of conducting 
constructions within a suitable topos. Logic is no longer something imposed 
from without; it becomes a ``physical variable''. A suspicion arises that at 
smaller and smaller distances, that we try to explore in our search for 
quantum gravity, this ``dependence on logic'' could be even more pronounced. 
The proposal has been put forward that at very short distances or high 
energies some regions of space-time are modeled by suitable toposes and 
specific quantum mechanical effects could be generated by their structural 
properties \cite{Krol2006}. 

In the present paper, we continue this line of research but our strategy is different. In category theory, logic and geometry cooperate with each 
other producing what is known under the name of Synthetic Differential 
Geometry (SDG). The main difference between SDG and standard differential 
geometry is that the  intuitionistic logic of SDG enforces the existence of 
infinitesimal objects which we symbolically denote by $D^n_{\infty}$ (for 
details see below). They essentially modify the local structure of spaces 
considered in SDG. This ``categorical geometry'' has already found several 
applications to general relativity \cite{Crane2006,Guts2001,MR-1991}. In the present paper, 
motivated by possible applications to the singularity problem and quantum 
gravity, we focus on an ``infinitesimal version'' of SDG, i.e. on SDG as it is 
restricted to $D^n_{\infty}$. We define (following \cite{SDG-2006,Lav1996,MR-1991}) an infinitesimal $n$-dimensional formal manifold and develop differential geometry on it. In particular, we show that the Riemann-Christoffel curvature tensor on locally $D^n_{\infty}$-spaces, as defined  by the Cartan's translation along the infinitesimal 2-chains \cite{MR-1991}, must be infinitesimal. 

To present SDG on infinitesimal spaces in a coherent way, a suitable topos must be specified which would provide a conceptual environment for the above interaction of logic, geometry 
and gravitation. Various choices are possible. Our choice falls on the topos 
${\cal G}$ of germ-determined ideals. It is a ``well adapted'' topos, in the 
sense that it ``smoothly'' generalizes all constructions required by the structure of space-time and general 
relativity. In particular, ${\cal G}$ contains a subcategory ${\mathbb{M}}$ of SET of smooth manifolds and their diffeomorphisms.

This change of perspective from the topos SET (the category with sets as 
objects and functions between sets as morphisms) to the topos $\cal{G}$ has 
far-reaching consequences both for cosmology -- by modifying the structure of 
space-time at small scales and possibly eliminating singularities, and for 
purely conceptual considerations -- it turns out to be related to the 
appearance of an exotic smooth structure on $\mathbb{R}^4$.

We proceed along the following line. In section 2, we prepare a terrain for 
dealing with geometry on an infinitesimal scale. In SDG one defines an 
$n$-dimensional formal manifold as a generalisation of the usual $C^{\infty 
}$-manifold. It is an object $M$ in a suitable category, in our case in the 
category $\cal{G}$, which has a cover $\{\varphi_i: U_i \rightarrow M\}$ by 
formally etal\'e monomorphisms, where $U_i$ are formal etal\'e subobjects in 
$\mathbb{R}^n$ (for details see below). The pairs $(U_i, \varphi_i)$ play the 
role of local charts on $M$. We need an infinitesimal version of this concept. 
It is elaborated in Definition 1.

In section 3, we start to develop differential geometry on an infinitesimal 
formal manifold. The main result of this section is the proof that on any 
infinitesimal formal manifold the curvature tensor assumes only infinitesimal values  on arbitrary 2-chains (Theorem 2). Something like that had to be expected, but the consequences of this theorem are unexpectedly 
far-reaching.

To unveil them, we construct, in section 4, a simplified quasi-cosmological 
model $S^3 \times \mathbb{R}$. We allow for $S^3$ to shrink to the zero size, 
thus producing a (topologically) cone singularity. This model can 
topologically be embedded in $\mathbb{R}^4$. We try to formally implement the 
following picture. We follow shrinking of $S^3$ to smaller and smaller sizes 
which causes the curvature of $S^3$ to grow dangerously. The process goes on, 
as usual, in the environment of the topos SET, but when the diameter of $S^3$ 
reaches a critical value $h$, the environment changes to that described by the 
topos $\cal{G}$, and `below $h$' everything is described in terms of 
$\cal{G}$. The curvature of $S^3$, instead of growing unboundedly, must now be 
infinitesimal, and the singularity, as it is expected in SET, is avoided in ${\cal G}$. 

Is this process of avoiding singularity totally invisible from the SET perspective? Not necessarily. It turns out that an observer in SET can effort a description that would take into account the existence of the critical value of $h$. However, such a description has to be done in a non-global way, i.e. by using two distinct coordinate patches, $\mathbb{R}^4_{<h}$ and $\mathbb{R}^4_{>h}$, belonging to a smooth atlas on $\mathbb{R}^4$. We call this type of evolution a hybrid evolution -- a hybrid since two toposes are engaged in it.

We meet here another unexpected effect. It is rather an elementary result that, given 
a smooth structure on $\mathbb{R}^4$, if there does not exist an open cover 
of it smoothly equivalent to the cover containing the single standard coordinate patch $\mathbb{R}^4$, this structure has 
to be exotic smooth (see e.g. \cite{Krol2016}). Consequently, the hybrid evolution, 
described above, has to be exotic smooth with respect to an exotic structure on 
$\mathbb{R}^4$, and this effect is due to an interaction between the infinitesimal and purely SET levels. In this way, we have a surprising result, namely making use of infinitesimals in $\cal{G}$ may have geometric consequences in SET.

The results obtained in this work suggest certain applications and 
provoke some comments which we include in section 5.

\section{Locally $D^n_{\infty}$-spacetimes}
\label{sec:1}
In this and subsequent sections we describe the construction of `locally 
infinitesimal' manifolds in a categorical setting. In particular, we are 
interested in the curvature of such manifolds. We follow the presentation of 
SDG given in \cite{SDG-2006} and \cite{MR-1991}. First, we introduce 
infinitesimals. Non-trivial infinitesimals do not exist in the category SET, but they do exist in 
some other toposes. As we remarked in the Introduction, we shall work with the 
topos $\cal{G}$ of germ-determined ideals in which infinitesimals do exist. It 
was introduced in \cite{Dubuc-1979} (see also \cite{Kock-1981}) and widely 
discussed in \cite{SDG-2006,MR-1991}. The importance of the topos ${\cal G}$ 
comes from the fact that it contains manifolds (i.e., there is an embedding of 
the manifold category into $\cal G$) and is closed under inverse image and 
exponentiations. As the consequence of the latter properties, $\cal G$ can 
contain spaces with various singularities and spaces of smooth functions. 
$\cal G$ is a Grothendieck topos defined on a site constructed of some reduced 
space of `smooth rings' $\mathbb{R}^n\to \mathbb{R}$ such that the $n$-potent 
infinitesimals, $n\in \mathbb{N}$, are modelled on the spectra of Weil 
algebras. Let us look at these concepts in some details.

Let $k$ be a commutative ring in SET. A Weil algebra $W=(k^n,\mu)$ is defined 
by the 2-linear mapping (multiplication) 
\begin{equation}\label{mu}
\mu : k^n\times k^n\to k^n 
\end{equation}
where $k^n$ is a commutative $k$-algebra with unit $(1,0,...,0)$, and there is 
the ideal $I\subset k^n$, defined by $\mu$, such that $I=(0,x_2,...,x_{n})$ 
and $I^n=0$ (this means that the $\mu$-product of $n$-elements of $k^n$ is 
$0$).  

Let further ${\cal E}$ be a Cartesian closed category and $R$ a commutative 
algebra over a $k$-object in ${\cal E}$ (the so-called $k$-algebra object in 
${\cal E}$). Basing on (\ref{mu}) one defines a commutative $R$-algebra 
$(R^n,\mu_R)$ in ${\cal E}$  with the unit $(1,0,...,0)  \in R^n$ and the 
$R$-algebra map $\pi:R^n\to R$ which is the projection on the first factor 
(this projection is called augmentation). The kernel of $\pi $, $\{x\in 
R^n|\pi(x)=0\}\simeq R^{n-1}$, is denoted by $I\otimes R$; it is the ideal in 
$ R\otimes W := (R^n,\mu_R)$. The $n$-fold powers of elements of $R\times I$, 
with respect to $\mu_R$, are equal to zero since $\mu_R$ is defined in terms 
of $\mu$, and both have the same `structure coefficients' \cite[p. 
62]{SDG-2006}. 

Turning now to SET, one can build \emph{finitely presented $k$-algebras $B$'s} 
which are given by the quotients
\begin{equation}\label{FT}
B=k[X_1,...,X_n]/(f_1(X_1,...,X_n),...,f_m(X_1,...,X_n))
\end{equation}
where $k[X_1,...,X_n]$ is the ring of $k$-coefficient polynomials in $n$ 
variables, and $(f_1(X_1,...,X_n),...,f_m(X_1,...,X_n))$ is the ideal spanned 
by the polynomials $f_i$. Again in the category ${\cal E}$, given the 
commutative $k$-algebra object $R$, one builds subobjects of $R\times ... 
\times R$ 
\begin{eqnarray*}
\lefteqn{Spec_RB =} \\
& & \{(r_1,...,r_n)\in R\times ... \times R = 
R^n|f_1(r_1,...,r_n)=0,...,f_m(r_1,...,r_n)=0\}.
\end{eqnarray*} 
If $FP\mathbb{T}_k$ is the category of finitely presented $k$-algebras then 
$Spec_RB$ can be regarded as a functor
\[ 
Spec_R: FP\mathbb{T}_k^{op}\to {\cal E} 
\] 
($op$ means, as usual, the opposite category in which morphisms are reversed) 
which to a $FP$ $k$-algebra $B$ assigns the object $Spec_RB$ such that 
$k[X]\to R$.\footnote{Notice that $k[X]$ is also a finitely presented 
$k$-algebra with the trivial ideal $(-)$ as a 0-generator.} Moreover, the 
functor $Spec_R$ preserves finite inverse limits \cite[p. 43]{SDG-2006}. 

It is well known that every Weil algebra $(k^n,\mu)$ over $k$ is finitely 
presented in terms of, say, $n$ generators\footnote{Any Weil algebra $W$ over 
$\mathbb{R}$ is a quotient $C^{\infty}(\mathbb{R}^n)/I_{k+1}$ where $I_{k+1}$ 
is the ideal generated by monomials of degree $k+1$ \cite[p. 160]{SDG-2006}.}. 
This is why presentation (\ref{FT}) applies also to them. For a $k$-algebra 
object $R$ in a category with finite inverse limits, like ${\cal E}$, which is 
Cartesian closed, one can build the objects $Spec_R(W)$. These objects are 
called \emph{infinitesimal objects} relative to $R$ in ${\cal 
E}$.\footnote{From now on we reserve the symbol $\mathbb{R}$ for the usual 
real line and use the symbol $R$ to denote the real line enriched by 
infinitesimals.} Given the projection-augmentation in every Weil algebra: 
$\pi: W\to k$ (i.e., the projection on the first factor of $k^n$ in 
(\ref{mu})), it can be seen that $Spec_Rk$ is a kind of terminal object 
$\textbf{1}$ amongst all Weil algebras in $FP{\mathbb{T}}_k$. Hence, in every 
infinitesimal object there is always a global element
\begin{equation}
\textbf{1}\to Spec_R(W).
\end{equation}

Let us consider a basic infinitesimal object $D\subset R$,  
$D=Spec_R(k[\epsilon])=\{r\in R|r^2=0\}$, where $k[\epsilon]=k[X]/(X^2)$ is
the Weil algebra with one generator, and $0: \textbf{1} \to D$ the 
canonical base point. There is no $0 \neq r \in \mathbb{R}$ in SET such that $r^2 =0$. Hence, $0: \mathbf{1} \rightarrow D$ is the only global element of $D$. Still, in $\cal{G}$ there are partial (i.e. non-global) elements corresponding to $d \in D, \, d \neq 0, \, d^2 =0$.

In general, any infinitesimal object $Spec_R(W)$ has a 
base point $b=Spec_R(\pi): Spec_R(k)\to Spec_R(W)$ where the functor $Spec_R$ 
is acting on the projection functor $\pi$ in the Weil algebra in 
$FP{\mathbb{T}}_k$. This can be presented in the diagram below (respecting the 
`$op$' order):
\begin{equation}
\xymatrix{
k\ar[d]_{{\pi}^{op}}_{}="0" && R= Spec_R(k)
           \ar[d]^{}_{}="1"\\
      \ar@{=>}"0";"1"^{{}_{b=SpecR(\pi)}}
W && Spec_R(W)
}
\end{equation}

Let us indicate some other infinitesimal objects that frequently appear in SDG 
constructions and are generated by spectra of Weil algebras
\begin{align}\label{eq2} 
\begin{split}
D_k=\llbracket x\in R|x^{k+1}=0 \rrbracket , \; k=1,2,3,...; \\ 
D_{k_1}\times D_{k_2}\times ... \times D_{k_n}\subset R^n \\
D(2)=\llbracket (x_1,x_2)\in R^2|x_1^2=x_2^2=x_1x_2=0 \rrbracket \\
D(n)=\llbracket (x_1,...,x_n)\in R^n|x_ix_j=0,\; \forall i,j=1,2,3,... 
\rrbracket \\
 D_k(n)=\llbracket (x_1,...,x_n)\in R^n|{\rm \, the\, product\, of \, any}\,  
k+1\, {\rm of} \, x_i\, {\rm is}\, 0 \, \rrbracket \\
 D_{\infty}^n=\bigcup_{k=1}^{\infty}D_k(n)\, .
\end{split}
\end{align}
These objects also contain partial elements in $\cal{G}$ besides the global one \cite{SDG-2006}.
The need for $D^n_{\infty}$ comes from the fact that $D$ is not an ideal in 
$R$ and from the following properties:
\begin{itemize}
\item[(a)] $D_{\infty}\subseteq R$ is an ideal (in the usual sense of ring 
theory).
\item[(b)] $D^n_{\infty}\subseteq R^n$ is a submodule.
\item[(c)] A map $t : D_{\infty}\to R$ with $t(0) = 0$ maps $D_k$ into $D_k$, 
for any k.
\end{itemize}

The following axiom expresses an important property of the infinitesimal 
objects $D_k(n)$ and maps $D_k(n)\to R$ in ${\cal G}$ (Axiom ($1''$), \cite[p. 
20]{SDG-2006}):
\begin{axiom}[1'']
For any $k=1,2,3,...$ and any $n=1,2,3,...$, every map $D_k(n) \to R$ is 
uniquely given by an $R$-polynomial in $n$ 
variables and of a degree not exceeding $k$.
\end{axiom}
The following corollary is the consequence of $(1'')$, assuming it holds for 
$R$ (Corollary 6.2, \cite[p. 20]{SDG-2006}, cf. (c) above):
\begin{corollary}\label{cor:1}
Every map $\phi: D_k(n)\to R^m$ with $\phi(0)=0$ factors through $D_k(m)$.
\end{corollary}

Let us now consider Weil algebra objects $ R\otimes W := (R^n,\mu_R)$ over $R$ 
in ${\cal E}$; $R$ is here a $k$-algebra object. To describe local   
$D^n_{\infty}$-manifolds in ${\cal E}$ we need yet another axiom of SDG to be 
satisfied by ${\cal E}={\cal G}$ \cite[p. 64]{SDG-2006}):
\begin{axiom}[$1^W_k$]\label{axiom:2}
For any Weil algebra $W$ over $k$, the $R$-algebra homomorphism $\alpha: 
R\otimes W\stackrel{\alpha}{\to } Spec_ R (W )$
is an isomorphism.
\end{axiom} 
Axiom ($1^W_k$) has important consequences. Let us notice that if $\alpha: 
R\otimes W\to  R^{Spec_ R (W )}$ is an isomorphism and $b 
=Spec_R(\pi):R=Spec_R(k)\to Spec_R(W)$ is sent by some $\beta \in 
R^{Spec_R(W)}$ to $0\in R$ (where $\pi$ is the projection in $W$) then, by the 
isomorphism $\alpha$, $\beta$ is sent to the ideal $R\otimes I$ in $R\otimes 
W$. This means that $\beta$ takes nilpotent values in $R$ for some $D_k$. 

Hence, we have:
\begin{proposition}\label{prop:1}[Proposition 16.3, p. 64]
Axiom ($1^W_k$) implies that any map
\[
Spec_R(W) \to R,  
\]
sending $b=Spec_R(\pi)$ to zero, factors through some $D_k$.
\end{proposition}
Now, we define local $D_{\infty}^n$-objects in ${\cal E}$ in analogy with 
external local $\mathbb{R}^n$-manifolds. We are interested 
in functorial properties of such objects, especially in their covering 
families of local `patches'. This is done with the help of Proposition \ref{prop:1} 
and Axiom \ref{axiom:2}. 

Let us consider a class ${\cal D}$ of morphisms in ${\cal G}$ which contain 
base points $\textbf{1}\to Spec_R(W)$ of infinitesimal objects $Spec_R(W)$. We 
define the class of ${\cal D}$-\'etale maps $f:M\to N$ in ${\cal G}$. Namely, 
$f:M\to N$ is ${\cal D}$-\'etale if for each $b: \textbf{1}\to Spec_R(W)$ the 
commutative square
\begin{equation}\label{diag2}
\xymatrix{
   M^{Spec_R(W)}
     \ar[rr]^{f^{Spec_R(W)}}
     \ar[dd]_{M^b}
      && N^{Spec_R(W)}
     \ar[dd]^{N^b} \\ \\
      M^{R}
     \ar[rr]^{f^R}
      && N^R}
\end{equation}
is a pullback square.
 
The following Theorem is crucial for the locally infinitesimal manifold concept:
\begin{theorem}
[\cite{SDG-2006}, Proposition 17.1 p. 70] \label{Th1}
The inclusion (a monic map) $(D_{\infty})^n \rightarrowtail R^n$ is a ${\cal 
D}$-{\'e}tale.
\end{theorem}

The structure of ${\cal D}$-{\'e}tale maps allows to consider $(D_{\infty})^n$ 
as formal $n$-dimensional objects, similarly to $R^n$-objects. This can be 
used to locally model more complicated `manifolds' in ${\cal E}={\cal G}$. 
Therefore,
\begin{definition}\label{def1}
An \emph{infinitesimal $n$-dimensional formal manifold} is an object $M$ in $\cal{E}$ for 
which there exists a `jointly epic' class of monic ${\cal D}$-\'etale maps
\[\{ (D_{\infty})^n_i \rightarrowtail M|i\in I\}. \]
\end{definition} 

A topology on such infinitesimal formal manifolds can be defined in terms of 
the relation $\sim_k$ on $M$ induced from $D_{k}(n)$ (see below footnote 5 and \cite[p. 74]{SDG-2006}). However, for our 
purposes there is no need to have any such topology on $M$. The possibility to 
work  with `manifolds' in ${\cal E}$ without defining any topology is a 
particular important feature of the categorical approach in SDG.

\section{The curvature of infinitesimal spacetimes}\label{curv}
\label{sec:2}
In this section, we describe the differential geometry on $n$-dimensional 
infinitesimal objects and determine their tensorial curvature. We assume that 
all objects are microlinear spaces in ${\cal G}$ (see the Appendix A). On the one hand, it should be expected that the `internal curvature' of any infinitesimal object should be infinitesimal if non-vanishing but, on the 
other hand, the tangent space at any point to infinitesimal formal 
$n$-manifold is still a vector space of dimension $n$ over $R$ (a module over $R$, see Appendix B).

Let $E\overset{p}{\to} M$ be a vector bundle on $M$. A connection $\nabla$ on 
$E$, $\nabla: M^D \times_M E\to E^D$ (see the Appendix B), determines the 
parallel transport of a vector $v\in E$ over a tangent vector by 
\[ \nabla(t,v)(0)=v,\; 0\in D, \] 
and the horizontal transport by
\[ p\circ \nabla(t,v)=t. \] 
The parallel transport should be linear in $v\in E$ and $t\in M^D$,
\[ 
\nabla(\alpha t,\beta v)(d)=\nabla(t,v)((\alpha\beta)\cdot d),\;\alpha, 
\beta\in R,\;d\in D.  
\] 
In this way, the transport $r_d(t,v)$ of $v$ along $t$ in an (infinitesimal) 
time period $d$ is defined. Let us notice that for the tangent bundle $E=M^D$ 
we have  $\nabla: M^D \times_M M^D\to M^{D\times D}$, so that
\begin{equation}\label{r}
r_d(t_1,v):= r_d(t_1,t_2)=\nabla(t_1,t_2)(d),\;t_1,t_2\in M^D,\; d\in D.
\end{equation}

Next, we generalize the transport over tangent vectors to the \emph{transport 
along the infinitesimal 2-chains}. Similarly as a tangent vector $t\in M^D$, 
the infinitesimal cell is given by the morphism $\gamma \in M^{D\times D}$. An 
infinitesimal 2-cube is thus (see Appendix B) \[ (\gamma, d_1,d_2)\in M^{D\times D}\times 
D\times D \] where $\gamma$ sends the infinitesimal ($d_1\times d_2$)-cube in 
$D^2$ into $M$.  In general, $n$-infinitesimal cube on $M$ is an element 
\begin{equation} 
(\gamma, d_1,...,d_n)\in M^{D^n}\times D^n. 
\end{equation}
The space of formal free modules over $R$ generated by infinitesimal 
$n$-chains is called the space of infinitesimal $n$-chains. 

In order to measure the curvature of the internal manifolds in ${\cal G}$, we 
define (following \cite{MR-1991}) a tensor $T(M)\times_MT(M)\times_MT(M)\to T(M)$. The value of this 
tensor is determined from a connection $\nabla$ on $M$. We should show that 
this tensor does not depend on $\gamma$, but only on $\nabla$.
The first step to do so is to define a map in ${\cal G}$ (\cite{MR-1991}, p. 
235)
\[
\overset{\approx}{R}:(M^{D^2}\times D \times D)\times_MM^D\to M^D 
\]
which to every 2-chain on $M$, $(\gamma,d_1,d_2)\in M^{D^2}\times D \times D$ 
and to a tangent vector $t_3\in M^D$ assigns the translated vector 
$\overset{\approx}{R}((\gamma,d_1,d_2),t_3)\in M^D$. The translation of $t_3$ 
is over the infinitesimal 2-chain $\gamma$. Thus given $\nabla$, parallel 
transport (\ref{r}) produces $r_{d_1}(t_1,t_2)(d_2)=\nabla(t_1,t_2)(d_1,d_2)$. 

Let us introduce an infinitesimal contour $\partial\gamma 
=\{\gamma_1,\gamma_2,\gamma_3,\gamma_4 \}$ around $D\times D$ centred at 
$(0,0)\in D^2$ by \[ \gamma_1=\gamma(-,0),\;\gamma_2=\gamma(d_1,-),\; 
\gamma_3=\gamma(-,d_2),\;\gamma_4=\gamma(0,-).\]
The parallel displacement of $t_3$ around $\partial\gamma$ gives
\begin{equation}
r(\gamma,d_1,d_2,t_3)=r^{-1}_{d_2}(\gamma_4,r^{-1}_{d_2}(\gamma_3,r_{d_2}(\gamma_2,r_{d_1}(\gamma_1,t_3)))),
\end{equation} 
so that the value of $\overset{\approx}{R}$ is just the difference of the 
translated vector and the initial vector $t_3$
\begin{equation}
\overset{\approx}{R}(\gamma,d_1,d_2,t_3)=r(\gamma,d_1,d_2,t_3)-t_3.
\end{equation} 
One then shows that $\overset{\approx}{R}(\gamma,d_1,d_2,t_3)=d_1 
d_2\overset{\sim}{R}$ where $\overset{\sim}{R}:M^{D^2}\times_MM^D\to M^d$ is 
the uniquely determined function. Namely, 
$\overset{\approx}{R}(\gamma,d_1,d_2,t_3)\in T_{\gamma(0,0)}(M)$ and, by 
fixing $\gamma$ with given $t_3$, we obtain 
$(\overset{\approx}{R})(d_1,d_2)\in T_{\gamma(0,0)}(M)$ i.e. 
$(\overset{\approx}{R}):D^2\to T_{\gamma(0,0)}(M)$. But $T_{\gamma(0,0)}(M)$ 
is microlinear, hence (\cite{MR-1991}, Proposition 1.4, p. 186)
 \[(T_{\gamma(0,0)}(M))^D\times 
(T_{\gamma(0,0)}(M))^D\overset{\alpha}{\to}(T_{\gamma(0,0)}(M))^{D\times D}  
\] is the isomorphism, and it gives the unique $\theta \in 
(T_{\gamma(0,0)}(M))^D$ such that $\theta(d_1\cdot 
d_2)=(\overset{\approx}{R})(d_1,d_2)\in T_{\gamma(0,0)}(M)$.  From the 
Kock-Lawvere axiom it follows that $\theta(d)=d\cdot v$ for the unique $v\in 
T_{\gamma(0,0)}$, hence
 \begin{equation}
 \theta(d_1\cdot d_2)=d_1\cdot d_2\cdot v= d_1\cdot d_2\cdot 
\overset{\sim}{R}(\gamma,t_3) \end{equation}
where $v=\overset{\sim}{R}(\gamma,t_3)$, and thus 
$\overset{\sim}{R}:M^{D\times D}\times_M M^D\to M^D$, as we have claimed.

\begin{definition}
[\cite{MR-1991}, p. 236]
The Riemann-Christofel tensor ${\cal R}: M^D\times_M M^D\times_M M^D\to M^D$ 
in ${\cal G}$ is given by the map:
\[{\cal R}(t_1,t_2)(t_3):=\overset{\sim}{R}(\nabla (t_1,t_2))(t_3). \] 
\end{definition}
This is well defined since $\nabla (t_1,t_2)\in M^{D\times D}$ and $t_3\in 
M^D$, so that ${\cal R}(t_1,t_2)(t_3)\in M^D$. One can show that 
$\overset{\sim}{R}(\gamma,t_3)$ does not depend on the 2-chain $\gamma$, but 
it does depend on the map $K:M^{D\times D}\to M^D\times_M M^D$ defining 
$\nabla$ (see Appendix B). 

The following lemmas lead to our main result.
\begin{lemma}
The non-zero curvature tensor ${\cal R}$ on the formal manifold $D_{\infty}$ 
assumes only infinitesimal values in some $D_k$. 
\end{lemma} 
\proof follows from the shape of the tangent space to $D_{\infty}$ which is 
$(D_{\infty})^D$, and thus ${\cal R}: (D_{\infty})^D\times_{D_{\infty}} 
(D_{\infty})^D\times_{D_{\infty}} (D_{\infty})^D\to (D_{\infty})^D$. This 
means that ${\cal R}(t_1,t_2)(t_3)\in (D_{\infty})^D$ and ${\cal 
R}(t_1,t_2)(t_3)(d)\in D_{k}(n)$ for some $k\in \mathbb{N}$ and $d\in D$. 
$\Box $

Let us recall some relations between infinitesimal objects (\ref{eq2})
\begin{align}\label{eq3}
\begin{split}
 D_{\infty}^n=\bigcup_{k=1}^{\infty}D_k(n),\; D_k(n)\subset D_l(n),k<l \\
 D_k(n) \subseteq (D_k)^n \\
 (D_k)^n \subseteq D_{ nk}(n)\,.
 \end{split}
 \end{align}

\begin{lemma}\label{l3}
The curvature tensor ${\cal R}$ on the formal manifold $D_{\infty}^n$ assumes 
only infinitesimal values in the object $D_k(n),n>1$ for some $k\in 
\mathbb{N}$. 
\end{lemma} 
\proof In this case, the tangent space is the object 
$(D^n_{\infty})^D=D^n_{\infty}\times D^n_{\infty}$ and thus ${\cal R}: 
(D^n_{\infty})^D\times_{D^n_{\infty}} (D^n_{\infty})^D\times_{D^n_{\infty}} 
(D^n_{\infty})^D\to (D^n_{\infty})^D$. Finally, from (\ref{eq3}) we have: 
${\cal R}(t_1,t_2)(t_3)\in (D^n_{\infty})^D$, and ${\cal 
R}(t_1,t_2)(t_3)(d)\in D_k(n)$ for some $k\in \mathbb{N}$, $d\in D$. $\Box $

Now, our main result concerning the curvature tensor of the locally 
$D^n_{\infty}$-formal manifolds (see Definition \ref{def1}):
\begin{theorem}\label{Th2}
The curvature tensor ${\cal R}$ of any locally $D^n_{\infty}$-formal manifold 
assumes only infinitesimal values in the object $D_k(m)$ for some $k\in 
\mathbb{N}$ and $m\geq n, m, n \in \mathbb{N}$.
\end{theorem}
\proof From Prposition \ref{prop:1} and Theorem \ref{Th1} it follows that, 
given $D^n_{\infty}$, it can be embedded in $R^n$ by a monic map and factor 
through some $D_k^n$ in $R^n$. Let us work with $D_k^n\subset R^n$ instead of 
$D^n_{\infty}$. In SET there is a Whitney embedding theorem for manifolds 
which states that any real smooth $n$-dimensional manifold $M^n$ can always be 
embedded in $\mathbb{R}^{2n}$. We claim that in ${\cal G}$, 
$D^n_{\infty}$-manifold, denoted by loc(n), is locally  monic-embeddable in 
$R^{2n}$ and it factors through $D^{2n}_l$ with some $l \geq k$. Let us consider 
$D^{2n}_{\infty}\rightarrowtail R^{2n}$, and notice that loc(n) in ${\cal G}$ 
is described by jointly epic family of (local) monic maps
\[\{ (D_{\infty})^n_i \rightarrowtail {\rm loc(n)}\;|i\in I\}. \] Given the 
embeddings $(D_{\infty})^n_i\rightarrowtail R^n_i,i\in I$, one builds a formal 
manifold $M^n$ by jointly epic family
\[ \{R^n_i\rightarrowtail M^n|i\in I  \}\,. \] 
However, $M^n$ in ${\cal G}$ can be obtained from SET with the help of the embedding functor 
$s: \mathbb{M} \rightarrow {\cal G}$ \cite{MR-1991}. Thus also the SET relation $M^n\subset 
\mathbb{R}^{2n}$ holds in ${\cal G}$ as $s(M^n)\subset R^{2n}$. From the 
construction it follows that ${\rm loc(n)}\subset D^{2n}_{\infty}$ which 
factors through some $D^{2n}_l, l \geq k$. However, from the Whitney theorem, 
the maximal dimension for the embedding is $2n$, which means that the 
infinitesimal space is such that $D^{m}_l, l \geq k, n \leq m\leq 2n$. From eqs. 
(\ref{eq3}) we deduce that the space $D_l(m)\subset D_l^m\subset D_{lm}(m),l\in \mathbb{N},2n \geq  m \geq n$ and $D_{\infty}^n=\bigcup_{k=1}^{\infty}D_k(n)$.
Since ${\rm loc(n)}\subset D^{2n}_{\infty}$, we can use Lemma $\ref{l3}$ to 
complete the proof. $\Box $

\section{A hybrid model}\label{s4}
In this section, we address the problem of SET-based constructions that would 
be sensitive to the existence of infinitesimal spaces in ${\cal G}$. The 
difficulty consists in the fact that the SET perspective causes unavoidable 
disappearance of non-zero infinitesimals. This is the consequence of the nonexistence 
of $D$-objects as subsets of $\mathbb{R}$ in SET. To overcome this difficulty 
we construct a hybrid model suitably combining both perspectives: the SET 
perspective and the $\mathcal{G}$ perspective.

Let $\overline{{B}^4}$ and $B^4$ be closed and open 4-balls (in 
$\mathbb{R}^4$), respectively. Then of course, $\partial \overline{{B}^4}\simeq 
S^3$. We consider a simplified model for an evolving universe given by (e.g. 
\cite{AsselKrol2014})
\[
S^3 \times \mathbb{R} 
\]
where, in analogy with the closed Friedman-Lema\^{i}tre cosmological model, 
$\mathbb{R}$ can be interpreted as a cosmic time and $S^3$ as a 3-dimensional 
instantaneous time section (although so far we remain on the purely 
topological level). The canonical relation holds
\[
S^3 \times \mathbb{R}\cup B^4= \mathbb{R}^4. 
\] 
$\mathbb{R}^4$  can be regarded as a Riemann manifold, and we can consider a 
smooth evolution in $\mathbb{R}^4$, 
\begin{equation}\label{evo1}
S^3 \times \mathbb{R}\subset \mathbb{R}^4, 
\end{equation} 
in the sense that the smooth evolution of $S^3\times \mathbb{R}$ respects the 
standard smoothness of $\mathbb{R}^4$.\footnote{The standard smooth structure of $\mathbb{R}^4$ is the unique structure in which the product $\mathbb{R} \times \mathbb{R} \times \mathbb{R} \times \mathbb{R}$ is a smooth product.} Now we allow for the smooth shrinking 
of the diameter $\rho_{S^3}$ of $S^3$ to the zero size (i.e., to the point 
${\rm pt.}\in \mathbb{R}^4$, which we situate at, say, $x_0=0, \, x_0 \in 
\mathbb{R}$).  Thus shrinking the size of $S^3$ to arbitrarily small values of 
the diameter $\rho_{S^3}$ is described as smooth contraction in the standard 
$\mathbb{R}^4$. Topologically, we have a cone over $S^3$ with the vertex ${\rm 
pt.}\in \mathbb{R}^4$. If we delete an open neighbourhood of the vertex, the 
cone becomes a standard smooth open 4-submanifold of $\mathbb{R}^4$ (without 
any `smoothing the corners' by isotopy). We call this vertex the singularity, 
but we should remember that it is a simple cone singularity rather than a 
curvature singularity met in standard cosmological models. Our aim is to 
prolong the evolution over this non-smooth vertex with the help of 
infinitesimally small elements. 

Switching between the categories SET and ${\cal G}$ (both of them are toposes) 
is in general governed by geometric morphisms that preserve much of the 
logical and intuitionistic set structures. However, there exists a special 
embedding of the category of smooth manifolds $\mathbb{M}$ into ${\cal G}$. 
Namely, we have
\begin{lemma}[\cite{MR-1991} Corollary 1.4, p.102]
The embedding of the category ${\mathbb{M}}$ into ${\cal G}$, 
$s:\mathbb{M}\hookrightarrow {\cal G}$, is full and faithful.  
\end{lemma}
In this way, one can do geometry `inside' $\cal{G}$ (more on the functor $s$ 
see below). In particular, the object $R_{\cal G}$ is $s(\mathbb{R})$, and 
similarly $N_{\cal G}=s(\mathbb{N})$. Also $s(M)^D\simeq s(TM)$ 
(\cite{MR-1991}, p. 111). However, not all `manifolds' that are internal in 
${\cal G}$, are an image of a manifold from SET by $s$. The important examples 
are infinitesimal spaces in ${\cal G}$ and locally $D^n_{\infty}$-formal 
manifolds. We want to find a SET-based manifestation of their existence. To 
this end we make the following assumptions.

Suppose that the continuous evolution (\ref{evo1}) is defined globally in SET, 
i.e. in the topological $\mathbb{R}^4$; moreover, 
\begin{quotation}
(A) there exists a scale $0\leq h\in \mathbb{R}$, below which (i.e. when the 
diameter $\rho_{S^3}<h$) the \emph{smooth} manifold $S^3\times \mathbb{R}$ 
is described internally in the topos ${\cal G}$, but `outside' the 4-ball 
$B^4$ (i.e. when the diameter $\rho_{S^3}>h$) $S^3 \times \mathbb{R}$ is 
the usual smooth manifold described in SET.

(B) for diameters $0\leq \rho_{B^4}<(h)$ smaller than some internal $(h)\in 
R_{{\cal G}}$, the internal \emph{smooth} manifold $S^3\times R_{\cal G}$ is a 
locally $D^4_{\infty}$-infinitesimal manifold (rather than the image 
$s(S^3\times \mathbb{R})$ under $s:\mathbb{M}\to {\cal G}$). 

(B') we do not decide what happens for negative $(h)$ (`from the other side of 
singularity'), whether the internal spheres $S^3$ are infinitesimal or not 
(this would depend on the particularities of a given model; anyway, the 
present model is only a toy model). 
\end{quotation}

These innocently looking assumptions have, in fact, dramatic consequences.

First, let us notice that, according to Theorem \ref{Th2},  the values of the 
Riemann tensor are infinitesimal on $D_k(m)$ for some $k\in \mathbb{N}$ and 
$m>4$. This means that when the contraction goes on, the 3-curvature of 
$S^3$ increases (as described in SET), and when the contraction crosses the 
scale $h \in \mathbb{R}$, it acquires its prolongation in $\mathcal{G}$ and, 
on the strength of assumption (B), `below' $(h)\in R_{\cal G}$ the components 
of the curvature become infinitesimal. Therefore, they assume values in the monad ${\cal M}_k\overset{\rm izo}{\simeq} D_k(m)$ \footnote{For any formal manifold
$M$ and any $x\in M$ , $M_k(x)$ is ``the k-monad around x'', i.e.  $M_k(x) := [[y \in M | x \sim_k y]]$ and $x \sim_k y \Leftrightarrow (x-y)\in D_k(m)$ for $R^m$.} rather than being arbitrarily large \cite[p.74]{SDG-2006}. In this way, the cone singularity has been avoided.

Second,  the existence of the `limiting values' $h\in \mathbb{R}$ and $(h)\in 
R_{\cal G}$ is paramount. If they exist then $h$ separates SET and ${\cal G}$ 
perspectives, whereas $(h)$ separates infinitesimal and non-infinitesimal 
descriptions. Reasoning exclusively in SET produces `singularity' (a violation of 
smoothness), whereas reasoning exclusively in ${\cal G}$ prevents to have a 
standard evolution. Therefore, both perspectives are indispensable.

What happens if one switches to SET but in such a way as to respect the existence of the 
separating $h\in \mathbb{R}$, and does this in a non-global way, i.e. below 
$h$ and above $h$ separately? The interaction between SET and $\mathcal{G}$ is 
governed by two functors, $s: \mathrm{SET} \rightarrow \mathcal{G}$ and  
$\Gamma: {\cal G}\to \mathrm{SET}$ (the latter is called global section). The 
global section functor $\Gamma $ has a left adjoint $\Delta$ which is the 
constant sheaf functor, so that $\Gamma \vdash \Delta$ and one has
\[
\xymatrix{
 {\mathbb{M}}\,\ar@<0.0ex>[r]^s  &  {\cal G}\, 
\ar@<-1.2ex>[r]^{\Gamma} & \, {\rm 
SET}\ar@<-1.2ex>[l]_{\Delta}.
}
\]
Moreover, $s: \mathbb{M}\to {\mathcal{G}}$ preserves transversal pullbacks, 
open covers, partitions of unity and compactness, and $\Gamma$, having left 
adjoint, preserves inverse limits, and it holds \cite[p. 226]{MR-1991}
\[\forall _{M\in \mathbb{M}}\; \Gamma(s(M))\overset{\rm diff}{\simeq} M,\; 
\forall_{M,N\in\mathbb{M}}\; \Gamma(s(N)^{s(M)})\overset{\rm diff}{\simeq} 
C^{\infty} (M,N).  \] 
Hence, $\Gamma$ and $s$ cancel each other but up to an isomorphism in 
$\mathbb{M}$ which is a \emph{smooth diffeomorphism} rather than just the 
identity diffeomorphism. This is in some sense a central datum for our 
construction.  It follows that there are two patches $\mathbb{R}^4_{<h}$ 
and $\mathbb{R}^4_{>h}$ in SET such that the first contains the image 
$\Gamma(S^3\times (-\infty,h)_{\cal G})$ in SET, and the second the image 
$\Gamma(s(S^3\times (h,\infty)))$. In this way, the separation `is visible' in 
SET.

The separation by $h$ means that one cannot have a single patch in 
$\mathbb{R}^4$ that would contain the whole of $S^3\times \mathbb{R}$. The 
question is whether one can glue these two patches in SET to obtain a  SET 
model for smooth evolution (\ref{evo1}), still taking into account the 
separation. First, let us consider the topological gluing. By slightly 
increasing both intervals in SET, i.e. $S^3\times (-\infty,h+\epsilon_1)$ (do not forget 
about (B')) and $S^3\times (h-\epsilon_2,\infty), \epsilon_{1,2}\in 
\mathbb{R}$, we obtain two coordinate patches in $\mathbb{R}^4$, each of them 
containing the increased subspaces. Let us use the same symbols as before for 
these patches, i.e. $\mathbb{R}^4_{<h}$ and $\mathbb{R}^4_{>h}$. In this 
way, we have two coordinate patches such that $\mathbb{R}^4_{<h}\cup 
\mathbb{R}^4_{>h}\overset{{\rm top.}}{\simeq}\mathbb{R}^4$. The topology of 
$\mathbb{R}^4$ is unique, therefore the topological gluing coincides with that 
of (\ref{evo1}). We now can consider the smooth evolution in SET, obtained 
from the topological evolution as described above, such that the patches 
$\mathbb{R}^4_{<h}$ and $\mathbb{R}^4_{>h}$ become standard smooth local 
coordinate patches of a smooth structure on $\mathbb{R}^4$.

A \emph{hybrid evolution}, or the evolution respecting the existence of separating $h\in \mathbb{R}$ and $(h)\in R_{\cal G}$, is thus a smooth evolution (\ref{evo1}) in $\mathbb{R}^4$ such that
\begin{itemize} 
\item[i.] the local coordinate patches $\mathbb{R}^4_{<h}$ and 
$\mathbb{R}^4_{>h}$ belong to some atlas of a smooth structure on 
$\mathbb{R}^4$, 
\item[ii.] no such atlas can be smoothly equivalent to the atlas with a single global smooth chart on 
$\mathbb{R}^4$.   
\end{itemize}

It is worth noticing that the hybrid evolution could also be defined with respect to a broader class of toposes. The candidate toposes are \emph{smooth toposes} which are models of SDG \cite{MR-1991}. Among them there exists the Basel topos \cite{MR-1991}. It was recently shown \cite{Krol2016} that the Basel topos indeed, by its very structure, realizes the hybrid geometry and can be used as a tool for distinguishing different, non-equivalent smooth atlases of coordinate patches on $\mathbb{R}^4$.

We can summarize the above in the form of the corollary:
\begin{corollary}
Smooth hybrid evolution in $\mathbb{R}^4$ is described with the 
help of an atlas on $\mathbb{R}^4$ necessarily containing at least two local coordinate patches.
\end{corollary}

It was for us a surprise to notice that this rather technically looking 
statement is related to deep results, obtained in 1980s, which revolutionized 
low dimensional geometry and topology. A remarkable theorem finalizing the effort of many mathematicians like Casson, Freedman, 
Donaldson, Taubes or Gompf (e.g. \cite{Gompf1985,StipGompf}), says that dimension 4 is distinguished from all 
other dimensions:
\begin{theorem}[Theorem 9.4.10 \cite{StipGompf}]
Only on $\mathbb{R}^4$ there exist uncountably many different pairwise 
nondiffeomorphic smooth structures. For any other $\mathbb{R}^n, n\neq 4$, there exists precisely one standard smoothness structure.
\end{theorem}
These 4-manifolds $\mathbb{R}^4$ that are not diffeomorphic to the standard 
smooth $\mathbb{R}^4$, but are all homeomorphic to it, are called exotic 
smooth $\mathbb{R}^4$. 
\begin{corollary}\label{cor3}
Each exotic $\mathbb{R}^4$ is a Riemannian smooth 4-manifold whose Riemann 
curvature tensor does not globally vanish.
\end{corollary}
\proof A smooth $\mathbb{R}^4$, on which Riemann tensor vanishes globally is 
flat and thus diffeomorphic to the standard $\mathbb{R}^4$. $\square$

Let us also quote an elementary but powerful lemma
\begin{lemma}[Lemma 9, \cite{Krol2016}] \label{Lem9}
Given a smooth structure on $\mathbb{R}^4$, if there does not exist any open
cover of $\mathbb{R}^4$ containing a single coordinate patch, 
this structure has to be exotic smooth.
\end{lemma}

With the above in mind the main results of this section can be formulated in 
the following way
\begin{theorem}\label{Th3}
For any $h\in \mathbb{R}$ and $(h)\in R_{\cal G}$ (as in assumptions (A) and 
(B)) such that the smooth evolution (\ref{evo1}) is hybrid, this evolution has 
to be exotic smooth with respect to some exotic $\mathbb{R}^4$.
\end{theorem}
\proof follows from Lemma {\ref{Lem9} and Corollary 2. $\Box $

We have here another instant of a subtle interaction between SET and 
$\mathcal{G}$, namely
\begin{theorem}\label{Riem}
The smooth hybrid evolution (with separating $h$ and $(h)$ as in Theorem 
\ref{Th3}) in dimension 4 gives rise to the nonvanishing Riemann tensor on a 
smooth $\mathbb{R}^4$. This Riemann tensor cannot be made zero by any 
diffeomorphism of $\mathbb{R}^4$.  
\end{theorem} 
\proof Theorem \ref{Th3} states that any hybrid evolution $S^3\times 
\mathbb{R}$ in SET with separating $h$ and $(h)$ has to be modeled on exotic 
$\mathbb{R}^4$. However, any exotic $\mathbb{R}^4$ cannot be globally flat, 
i.e. the Riemann tensor cannot globally  vanish (Corollary \ref{cor3}). 
Otherwise there would be a diffeomorphism of exotic and standard 
$\mathbb{R}^4$. $\square$ 

We have here a conceptually interesting result: the modification of smoothness 
on $\mathbb{R}^4$ (in SET perspective) is driven by what happens on the 
infinitesimal level (in $\mathcal{G}$ perspective).  

\section{Applications and comments}
In this paper, we have focused on differential geometry on the smallest 
possible -- infinitesimal -- scale. Besides of being interesting in itself, 
our results could naturally be expected to have important applications. As far 
as physical applications are concerned two of them seem to be especially 
attractive -- the singularity problem in cosmology and the problem of gravity 
on the Planck level. In both these problems curvature of space-time is 
involved. 

If, in contrast with what is predicted by general relativity, on approaching 
singularity the curvature, instead of unboundedly growing, becomes infinitesimal, 
then even the strongest singularities can be avoided. Moreover, since in SDG 
every function is smooth, one could smoothly join a contracting phase of the 
universe, through the almost-singularity, to its expanding phase. It is true that 
our results concern only Riemannian manifolds and in cosmology one deals with 
pseudo-Riemannian manifolds, but one could expect that also in this case strict 
results will be analogous. Moreover, one should take into account that the 
infinitesimal objects have no points individualised by their usual real 
coordinates, and the problem of the metric structure could demand careful 
rethinking. It might be that on this level there is simply no distinction 
between space and time directions (as it is suggested by some approaches to 
quantum gravity).

It is tempting to identify the infinitesimal level with the level beyond the 
Planck threshold (this suggestion is intimated in choosing the letter $h$ to 
denote the scale distinguishing the regimes SET and $\mathcal{G}$), but at this 
stage of investigation it is certainly premature. In the present work, our 
hybrid model plays rather a heuristic role. However, it has recently been shown that a cosmological model $S^3\times \mathbb{R}$ with an exotic smooth structure generates realistic parameters for the cosmological inflation (e.g. \cite{TAMJK2013,AsselKrol2014}).  Some other interesting results \cite{TAMJK2013,Freedman1979} state that exotic $\mathbb{R}^4$ and $S^3\times \mathbb{R}$ lead to the wild embedding $S^3\overset{\rm wild}{\hookrightarrow} \mathbb{R}^4$ (similarly to the well known embedding of the Alexander sphere $S^2\overset{\rm wild}{\hookrightarrow} \mathbb{R}^4$ \cite{Alex1924}) and to quantum non-commutative $C^{\ast}$-algebra \cite{TAMJK2013,Turaev1991}. In the light of these result, the fact that our hybrid model has revealed a connection with a smooth exotic structure could also be significant. It would be interesting to study how does the change from a tamed embedding to a wild embedding influence the appearance of infinitesimals and the `below $h$' scale. These topics certainly open a promising area for the future research.

Finally, it should be stressed that all the above results depend on changing the role of logic in our doing science: instead of being an  \textit{a priori} judge of our theories, logic changes into a `physical variable'. It seems that we now are at the threshold of another conceptual revolution.

\section*{Appendix A}
Many arguments in SDG are based on the microlinear property of formal 
manifolds. In this Appendix, we briefly introduce this concept. Our presentation follows that of \cite{SDG-2006,MR-1991}. Roughly speaking, microlinear space is a space that behaves, with respect to maps from infinitesimal spaces to itself, as it had local coordinates. The following considerations lead to the precise definition.

For any Cartesian closed category ${\cal E}$, the functor $F_M:{\cal E}\to 
{\cal E}$, defined for an object $M\in {\cal E}$ by $X\to M^X$, is the 
contravariant functor that sends colimit diagrams in ${\cal E}$ into the limit 
diagrams, i.e. 
\[ M^{{\displaystyle{\lim_{\overset{\longrightarrow}{i}}}}X_i}\simeq 
\lim_{\overset{{\longleftarrow}}{i}}X_i. \]
It can happen that $X_i\to X$ is not a colimit diagram but $M^X\to M^{X_i}$ is 
the limit cone. If this is not the case for infinitesimal objects $X_i$ in ${\cal 
E}$, then the object $M$ is said to be a microlinear object. More precisely, 
let us consider all limit diagrams in the category of commutative 
$\mathbb{R}$-algebras whose vertices are Weil algebras. Let us apply 
$Spec_R$-functor to obtain the class $Cocone(Spec_R(Weil))$ of cocones from 
the above class of limit Weil-algebra diagrams. They are not necessarily colimit 
diagrams. 
\begin{definition}
The object $M\in {\cal E}$ is \emph{microlinear} if every cocone from \\ 
$Cocone(Spec_R(Weil))$ becomes a colimit diagram under the action of the 
functor $M^{(-)}$.  
\end{definition}
Our basic example is given by the pullback diagram of Weil algebras
\begin{equation}\label{diagA1}
\xymatrix{
   \mathbb{R} && \mathbb{R}[\epsilon_2]
      \ar[ll] \\ \\
      \mathbb{R}[\epsilon_1]
     \ar[uu] && \mathbb{R}[\epsilon_1,\epsilon_2]
      \ar[ll]  \ar[uu]
     }
\end{equation}
which, by the functor $Spec_R$, is sent to
\begin{equation}\label{diagA2}
\xymatrix{
   \textbf{1}=Spec_R(\mathbb{R})
     \ar[rr]
     \ar[dd]
      && D=Spec_R(\mathbb{R}[\epsilon_1])
     \ar[dd]  \\ \\
      D=Spec_R(\mathbb{R}[\epsilon_2])
     \ar[rr]
      && D(2)=Spec_R(\mathbb{R}[\epsilon_1,\epsilon_2])}
\end{equation}
It is not a pushout diagram. However, the functor $R^{(-)}$ takes it into 
the pullback diagram
\begin{equation}\label{diagA3}
\xymatrix{
   R && R\times R
      \ar[ll] \\ \\
      R\times R
     \ar[uu] && R\times R\times R
      \ar[ll]  \ar[uu]
     }.
\end{equation}
This follows from the Kock-Lawvere axiom since then $R^D=R\times R$ and 
$R^{D\times D}=R\times R\times R$.

Let us suppose that $k$ is a field; then from Axiom ($1^W_k$) we have in 
${\cal G}$:
\begin{proposition}[Proposition D.1, \cite{SDG-2006}]
$R_{\cal G}$ in ${\cal G}$ is microlinear.
\end{proposition}
Microlinear spaces comprise a full subcategory in ${\cal G}$ which is closed 
under inverse limits and exponentiation. Moreover, all manifolds in SET, by 
the embedding functor $s: \mathbb{M}\to {\cal G}$, i.e. by $s(M)$, are 
microlinear in ${\cal G}$ (\cite{MR-1991}, p. 227). 
Infinitesimal spaces in ${\cal G}$ which are not images under $s$ are also microlinear
\begin{proposition}
$D_{\infty}^n\subset R_{\cal G}^n,n=1,2,...$ are microlinear in ${\cal G}$.
\end{proposition}
\proof The result follows from Proposition 17.6 \cite[p.74]{SDG-2006} and Theorem \ref{Th1}. $\square$

\section*{Appendix B}
In this Appendix, we briefly present some elements of intuitionistic differential geometry on formal 
$n$-manifolds (which are locally $(D_{\infty})^n$-spaces) in the category ${\cal G}$. They are 
important for deriving results of Sec. \ref{curv}.  We  mostly follow 
\cite{Lav1996}.  Let $i:D(2)\to D\times D$ be the canonical injection of the 
infinitesimal spaces, and let $M$ be a microlinear object in ${\cal G}$, then one 
has the canonical mapping $M^i: M^{D\times D}\to M^{D(2)}$. This mapping and the 
isomorphism $M^{D(2)}\simeq M^D\times_M M^D$ determine the morphism 
\begin{equation}\label{eq8} 
K: M^{D\times D}\to M^D\times_M M^D.
\end{equation} 
The connection $\nabla:M^D\times_M M^D\to 
M^{D\times D}$ on a microlinear space $M$ is the section of $K$ with the 
following linearity conditions holding for every $(t_1,t_2)\in M^D\times_M 
M^D$ and $d_1,d_2\in D, \alpha\in R_{\cal G}$
\begin{align}\label{eq5}
\begin{split}
\nabla(t_1,t_2)(d_1,0)=t_1(d_1), \\
 \nabla(t_1,t_2)(0,d_2)=t_2(d_2), \\
 \nabla(\alpha\cdot t_1,t_2)(d_1,d_2)=\nabla(t_1,t_2)(\alpha\cdot d_1,d_2),\\
 \nabla(t_1,\alpha\cdot t_2)(d_1,d_2)=\nabla(t_1,t_2)(d_1,\alpha\cdot d_2).
 \end{split}
 \end{align} 
A vector tangent to a microlinear $M$ at $m\in M$ is $t\in M^D$ such that $t(0)=m$. The space of all tangent vectors $T_mM$ to $M$ at $m\in M$ constitutes a module over $R$ \cite[Proposition 1, p. 62]{Lav1996}. 
A microsquare $\gamma: D\times D\to M$ is the tangent to $M^D$ since $M^{D\times D}\simeq (M^D)^D$ is given by $\tau(d_1)(d_2)=\gamma(d_1,d_2)$, where $\tau \in (M^D)^D$.
For $K$ as in (\ref{eq8}), we have $K:\gamma \to (t_1,t_2)$ and $t_1(d)=\gamma(d,0),\; t_2(d)=\gamma(0,d)$.

The infinitesimal parallel transport associated with the connection $\nabla$ is now defined by the mappings \begin{align*}
\begin{split}p_{(t,e)}(t')(d)=\nabla(t,t')(e,d) \\ 
q_{(t,e)}(t_1)(d)=\nabla (\nabla(t,t)(e,\cdot ),t_1)(-e,d) \end{split}
 \end{align*}  
where $t\in T_{t(0)}M, e, d\in D$ and $t_1\in T_{t(e)}M$ so that 
\begin{align*}
\begin{split} p(t,e): T_{t(0)}M\to T_{t(e)}M\\
 q(t,e): T_{t(m)}M\to T_{t(0)}M.\end{split}
 \end{align*} 
\begin{proposition}[\cite{Lav1996}, Proposition 6, p. 164]
Mappings $p(t,e)$ and $q(t,e)$ are inverse isomorphisms between fibers $T_{t(0)}M$ and $T_{t(e)}M$.
\end{proposition}



\end{document}